\documentclass[12pt,reqno]{amsart}
\newtheorem{theorem}{Theorem}
\newtheorem{proposition}[theorem]{Proposition}
\newtheorem{corollary}[theorem]{Corollary}
\newtheorem{remark}{Remark}

\newfont{\bb}{msbm10 at 12pt}

\def\pf{{\textit {Proof :} }}
\def\qed{\hfill{q.e.d.}\smallskip\smallskip}

\def\Si{{\Sigma}}     
\def\ga{{\gamma}}     
     
\def\Ga{{\Gamma}}

\def\phi{{\varphi}}

\def\TM{{\mathrm{T}\mathrm{M}}}  
   
\def\M{\mathrm{M}}
\def\T*M{\mathrm{T}^*\mathrm{M}}
\def\TBM{\mathrm{T}(\pa\M)}

\def\gna{\na^{\mathbf{S}}}
\def\DS{\D^{\mathbf{S}}}
\def\H{\mathrm{H}}
\def\R{\mathrm{R}}

\def\bg{\bar{g}}
\def\L{\mathrm{L}}
\def\B{\mathrm{B}}
\def\MIT{\mathrm{MIT}}
\def\CHI{\mathrm{CHI}}
\def\L{\mathrm{L}}
\def\BMIT{\mathbb{B}_{\MIT}}
\def\BCHI{\mathbb{B}_{\CHI}}
\def\Pe{\mathrm{P}}

\def\D{\mathrm{D}}
\def\hs{\mathbb{S}^n_+}
\def\S{\mathbf{S}}
\def\DSC{{\overline{\D}}^{\S}}
\def\V{\mathrm{V}}
\def\G{\mathrm{G}}
\def\ov{\overline}
\def\F{\mathrm{F}}

\let\pa\partial     
\let\na\nabla     
     
\DeclareMathAlphabet{\doba}{U}{msb}{m}{n}

\def\End{{\mathop{\rm End}}}     
\def\Hom{{\mathop{\rm Hom}}}     
\def\Spin{{\mathop{\rm Spin}}}     
\def\SO{{\mathop{\rm SO}}}   
\def\dim{{\mathop{\rm dim}}}  
     
\def\ker{{\mathop{\rm Ker}}} 
\def\im{{\mathop{\rm Im}}}

\def\aire{\mathrm{Area}(\M^2,g)}
     
\def\Id{{\mathop{\rm Id}}}

\let\<\langle     
\let\>\rangle

\long\def\komment#1{} 

\begin{document}


\title{The Hijazi Inequality on Manifolds with Boundary}     
\author{Simon Raulot}     
\date{\today}


\maketitle

\begin{abstract}
In this paper, we prove the Hijazi inequality on compact Riemannian spin
manifolds under two boundary conditions: the condition associated with a chirality
operator and the Riemannian version of the $\MIT$ bag condition. We then
show that the limiting-case is characterized as being a half-sphere for the
first condition whereas the equality cannot be achieved for the second.
\end{abstract}


\section{Introduction}


In \cite{friedrich}, T. Friedrich established an inequality relating
the eigenvalues of the Dirac operator on a compact $n$-dimensional
Riemannian spin manifold without boundary with its scalar curvature $\R$. This
inequality is given by
\begin{eqnarray}
\lambda^2\geq\frac{n}{4(n-1)}\,\R_0,
\end{eqnarray}

where $\R_0$ denotes the infimum of the scalar curvature of $\M$. 
For $n\geq 3$, using the conformal covariance of the Dirac operator O. Hijazi \cite{hijazi:86} proved that
\begin{eqnarray} 
\lambda^2\geq\frac{n}{4(n-1)}\,\mu_1(\L), 
\end{eqnarray}

where $\mu_1(\L)$ is the first eigenvalue of the conformal Laplacian given by
$\L=\frac{4(n-1)}{n-2}\Delta+R$. The operator $\L$ is a second order
conformally covariant differential operator relating the scalar curvatures of
two metrics in the same conformal class. In \cite{hijazi:91}, O. Hijazi derives a
conformal lower bound for any eigenvalues $\lambda$ of the Dirac
operator involving the Yamabe invariant
$\mu(\M)$. Indeed, he proved that if $n\geq 3$, then
\begin{eqnarray}\label{hjc}
\lambda^2\,\mathrm{vol}(\M,g)^{\frac{2}{n}}\geq\frac{n}{4(n-1)}\,\mu(\M).
\end{eqnarray}

The Yamabe invariant $\mu(\M)$ has been introduced 
in \cite{yamabe:60} in order to solve the following problem 
now called the {\it Yamabe problem:} given a closed compact Riemannian
manifold $(\M^n,g)$, is there a metric in the conformal class of $g$ such
that the scalar curvature is constant? For $n=2$, C. B\"ar \cite{baer:92} showed that 
\begin{eqnarray}
\lambda^2\geq\frac{2\pi\,\chi(\M^2)}{\mathrm{Area}(\M^2,g)},
\end{eqnarray} 

where $\chi(\M^2)$ is the Euler characteristic class. A natural question is
then to ask if those results still hold if we
consider manifolds with boundary. In \cite{hijazi.montiel.roldan:01}, the authors prove Friedrich-type inequalities under four elliptic boundary conditions and under some curvature assumptions
(the non-negativity of the mean curvature). Two of these boundary
conditions are (global) Atiyah-Patodi-Singer type conditions and two local
boundary conditions. The present paper being devoted to the conformal
aspect of those results, the choice of the boundary condition is 
important. As pointed out in \cite{hijazi.montiel.zhang:2}, the
Atiyah-Patodi-Singer type conditions are not conformally invariant, while
the local boundary conditions are indeed conformally invariant. 
Moreover we don't assume that the boundary $\pa\M$ has nonnegative
mean curvature and then we prove:

\begin{theorem}\label{mainth}
Let $(\M^n,g)$ be an $n$-dimensional connected compact Riemannian spin manifold with
non-empty boundary $\pa\M$. For $n\geq 3$ if $\M$ has positive Yamabe invariant
and for $n=2$ if $\M$ is a surface of genus $0$ whith compact connected
boundary, then:
\begin{itemize}

\item[ 1)] Under the $\CHI$ boundary condition associated with a chirality
  operator (see Section~\ref{chibound}), the
  spectrum of the Dirac operator $\D$ of $\M$ is a sequence of
  unbounded real numbers $\{\lambda^{\CHI}_k\;/\;k\in\mathbb{Z}\}$. For
  $n\geq 3$, any eigenvalue $\lambda^{\CHI}$ of the Dirac operator satisfies
\begin{eqnarray*}
\big(\lambda^{\CHI}\big)^2\geq\frac{n}{4(n-1)}\,\mu_1(\L).
\end{eqnarray*}

\noindent For $n=2$, we have
\begin{eqnarray*}
\big(\lambda^{\CHI}\big)^2\geq\frac{2\pi}{\aire}.
\end{eqnarray*}

\noindent Moreover, equality holds if and only if the manifold $\M$ is isometric to the half-sphere $\hs(r)$ with radius $r$, where $r$ depends on
the first eigenvalue of the Dirac operator under this boundary condition.

\item[ 2)] Under the $\MIT$ boundary condition (see Remark~\ref{rem2}), the
spectrum of the Dirac operator $\D$ is an unbounded
discrete set of complex numbers $\{\lambda^{\MIT}_k\;/\;k\in\mathbb{Z}\}$
with positive imaginary part. For $n\geq 3$, any eigenvalue $\lambda^{\MIT}$ of the Dirac
operator satisfies
\begin{eqnarray*}
|\lambda^{\MIT}|^2>\frac{n}{4(n-1)}\,\mu_1(\L). 
\end{eqnarray*}

\noindent For $n=2$, we have 
\begin{eqnarray*}
|\lambda^{\MIT}|^2>\frac{2\pi}{\aire}.
\end{eqnarray*}
\end{itemize}

\noindent The real number $\mu_1(\L)$ is the first eigenvalue of the eigenvalue
boundary problem 
$$\left\lbrace
\begin{array}{ll}
\L u=\mu_1(\L)u & \quad\mathrm{on}\;\M\\
\B u=0 & \quad\mathrm{along}\;\pa\M,
\end{array}
\right.$$

with $\B$ the mean curvature operator acting on
smooth fonctions on the manifold $\M$. 
\end{theorem}

Finally, we extend Inequality (\ref{hjc}) to the case of manifolds with boundary. The author would like to thank the referee for helpful comments.


\section{Spin Manifolds with Boundary}


In this section, we summarize some basic facts about spin manifolds with
boundary. Standard references on this subject can be found in
\cite{hijazi.montiel.zhang:2}. Let $(\M^n,g)$ be an $n$-dimensional
Riemannian spin manifold with boundary and denote by $\na$ the Levi-Civita
connection on the tangent bundle $\TM$ and the associated bundles. We choose a spin structure and denote by
$\Spin(\M)$ the principal bundle with structural group $\Spin_n$ given by
this spin structure. The spinor bundle on the manifold $\M$ is then the
complex vector bundle of rank $2^{[\frac{n}{2}]}$, denoted by $\Sigma\M$,
associated with the complex spinor representation. This representation provides a Clifford multiplication
$$\gamma :\mathbb{C}l\,(\M)\longrightarrow\End(\Sigma\M),$$

which is a fibre preserving algebra morphism. The spinor bundle $\Sigma\M$
is endowed with a natural hermitian product, denoted by $\<\;,\;\>$, and with a
spinorial Levi-Civita connection $\na$ acting on spinor fields,
i.e. on sections of the spinor bundle (see \cite{lawson.michelson:89} or \cite{fried}, for
example). We can easily show that the spinorial connection
$\na$ is compatible with the hermitian product $\<\;,\;\>$, i.e.
\begin{eqnarray}\label{sl}
X\<\psi,\varphi\>=\<\na_X\psi,\varphi\>+\<\psi,\na_X\varphi\>,
\end{eqnarray}
 
for all $\psi,\varphi\in\Ga(\Sigma\M)$ and for all
 $X\in\Ga(\TM)$. Moreover, they also satisfy the following properties:
\begin{eqnarray}
\label{skew1}
\<\gamma(X)\psi,\gamma(X)\varphi\>=g(X,X)\<\psi,\varphi\>\\
\label{skew2}
\na_X\big(\gamma(Y)\psi\big)=\gamma(\na_X Y)\psi+\gamma(Y)\na_X\psi,
\end{eqnarray}
 
for all $\psi,\varphi\in\Ga(\Sigma\M)$ and for all $X,Y\in\Ga(\TM)$. 
Note that Identity (\ref{skew1}) implies that Clifford multiplication by a unit tangent vector field is skew-symmetric on $\Sigma\M$. The Dirac operator
$\D$ on $\Sigma\M$ is then the first order elliptic operator locally given by
$$\begin{array}{lcll}
\D: & \Ga(\Sigma\M) & \longrightarrow & \Ga(\Sigma\M)\\
& \psi & \longmapsto & \sum_{i=1}^{n}\gamma(e_i)\na_{e_i}\psi,
\end{array}$$

where $\{e_1,...,e_n\}$ is a local orthonormal frame of $\TM$. 

\noindent Consider now the boundary $\pa\M$ which is an oriented Riemannian
hypersurface of $\M$ with induced orientation and Riemannian
structure. Then there exists a unit vector field $\nu$ normal to the
boundary which allows to pull-back the spin structure over $\M$ to a spin
structure over the boundary $\pa\M$. Hence we have that the restriction
\begin{eqnarray*}
\mathbf{S}(\pa\M):=\Sigma\M_{|\pa\M}
\end{eqnarray*}

\noindent is a left module over $\mathbb{C}l\,(\pa\M)$ with Clifford multiplication
\begin{eqnarray*}
\ga^{\S}:\mathbb{C}l\,(\pa\M)\longrightarrow\End(\S(\pa\M))
\end{eqnarray*}

\noindent given by $\ga^{\S}(X)\psi=\ga(X)\ga(\nu)\psi$ for all
$X\in\Gamma(\TM)$ and $\psi\in\Gamma\left(\S\left(\pa\M\right)\right)$. Now
let $\na^{\pa\M}$ be the Levi-Civita connection of the boundary $(\pa\M,g)$
 and let $(e_1,...,e_{n-1},e_n=\nu)$ be a local orthonormal frame of $\TM$, then the
Riemannian Gauss formula states that for $1\leq i,j\leq n-1$,
\begin{eqnarray*}
\na_{e_i}e_j=\nabla_{e_i}^{\pa\M}e_j+g(Ae_i,e_j)\nu,
\end{eqnarray*}

where $AX=-\na_X\nu$ is the shape operator of the boundary $\pa\M$. We can
then relate the two associated spin connections. Indeed, if $\na$ (resp. $\gna$) is the Levi-Civita connection on the spinor bundle $\Sigma\M$
(resp. $\bf{S}(\pa\M)$), we have the spinorial Gauss Formula (for more
details, see \cite{trautman}, \cite{baer:98} or \cite{morel}):
\begin{eqnarray*}
(\na_X\psi)_{|\pa\M}=\gna_X\psi_{|\pa\M}+\frac{1}{2}\gamma^{\S}(AX)\psi_{|\pa\M}.
\end{eqnarray*}

for all $X\in\Ga\big(\TBM\big)$ and for all $\psi\in\Ga(\Sigma\M)$. The spinor bundle
$\mathbf{S}(\pa\M)$ is also endowed with a
Hermitian metric denoted by $\<\;,\;\>$ induced from that on
$\Sigma\M$. The induced metric, the Clifford multiplication and  the Levi-Civita
connection satisfy properties (\ref{sl}), (\ref{skew1}) and (\ref{skew2}),
i.e. the spinor bundle $\S(\pa\M)$ is a Dirac bundle. 

\noindent The induced spin structure on the boundary allows to construct an
intrinsic spinor bundle $\Sigma(\pa\M)$ over $\pa\M$. This bundle is naturally endowed with
a Hermitian metric, a Clifford multiplication $\ga^{\pa\M}$ and a spinorial
Levi-Civita connection $\na^{\pa\M}$. It is not difficult to show that this
bundle can be identified with the restricted spinor bundle $\S(\pa\M)$ if
$n$ is odd. In this case, the Clifford multiplication and the Levi-Civita 
connection on $\Sigma(\pa\M)$ correspond to
the Clifford multiplication and the Levi-Civita connection on $\S(\pa\M)$. If $n$
is even, the spinor bundle $\S(\pa\M)$ could be identified with the direct
sum $\Sigma(\pa\M)\oplus\Sigma(\pa\M)$. Moreover, the Clifford multiplication
$\ga^{\S}$ correspond to $\ga^{\pa\M}\oplus-\ga^{\pa\M}$ and the
Levi-Civita connection $\na^{\S}$ to $\na^{\pa\M}\oplus\na^{\pa\M}$.

\noindent We can now define several Dirac operators acting on sections of these
bundles (for a complete review of these operators, see \cite{morel}). However, in our case, the most
important one is the {\it{boundary Dirac operator}} acting on sections of
$\S(\pa\M)$. This operator is given by composition of the Clifford
multiplication $\ga^{\S}$ and the spinorial Levi-Civita connection $\gna$. This
operator is denoted by $\DS$ and is locally given by
\begin{eqnarray*}
\DS(\psi)=\sum_{i=1}^{n-1}\ga^{\S}(e_i)\gna_{e_i}\psi.
\end{eqnarray*}

\noindent Note that in the case of a closed compact manifold
without boundary, the classical Dirac operator has exactly one self-adjoint $\mathrm{L}^2$ extention,
so it has real discrete spectrum. In the case of a manifold with boundary, a defect of symmetry appears, given by the formula 
\begin{eqnarray}\label{ipp}
\int_{\M}\<\D\psi,\varphi\>dv(g)-\int_{\M}\<\psi,\D\varphi\>dv(g)=-\int_{\pa\M}\<\gamma(\nu)\psi,\varphi\>ds(g),
\end{eqnarray}

for $\psi,\varphi\in\Ga(\Sigma\M)$, and where $\nu$ is the inner
unit vector field along the boundary and $dv(g)$ (resp. $ds(g)$) is the
volume form of the manifold $\M$ (resp. the boundary $\pa\M$). According to
this formula, we note that the Dirac operator $\D$ is not symmetric,
but we will see that, under suitable boundary conditions, the l.h.s of
(\ref{ipp}) vanishes.

\noindent In order to prove Theorem~\ref{mainth} which is an estimate of the
fundamental Dirac operator eigenvalues of the ambient manifold $\M$ under 
suitable boundary conditions, we will give an inequality called
``spinorial Reilly inequality'' relating the geometry of the manifold $\M$
and that of its boundary $\pa\M$ (see \cite{hijazi.montiel.roldan:01}, 
\cite{hijazi.montiel.zhang:1} or \cite{hijazi.montiel.zhang:2}). The spinorial Reilly inequality
is based on the Schr\"odinger-Lichnerowicz formula, given by
\begin{eqnarray*}
\D^2=\na^*\na+\frac{1}{4}\R,
\end{eqnarray*}

\noindent where $\R$ is the scalar curvature of $\M$. An integral version of this formula leads to (see
\cite{hijazi.montiel.roldan:01} for a proof of the following proposition):

\begin{proposition}\label{p4}
For all spinor fields $\psi\in\Ga(\Sigma\M)$, we have:
\begin{eqnarray}\label{ir}
\int_{\pa\M}\big(\<\DS\psi,\psi\>-\frac{n-1}{2}\H|\psi|^2\big)ds(g)\geq\int_{\M}\big(\frac{1}{4}\R|\psi|^2-\frac{n-1}{n}|\D\psi|^2\big)dv(g),
\end{eqnarray}

where $\R$ is the scalar curvature of the manifold $\M$,
$\H=\frac{1}{n-1}\mathrm{tr}(A)$ is the mean curvature of the boundary and
$dv(g)$ (resp. $ds(g)$) is the Riemannian volume form of $\M$ (resp. $\pa\M$). Moreover
equality occurs if the spinor field $\psi$ is a twistor-spinor, i.e. if it
satisfies $\Pe\psi=0$ where $\Pe$ is the twistor operator
acting on $\Sigma\M$ which is locally given for all $X\in\Ga(\TM)$ by:
\begin{eqnarray}
\Pe_X\psi=\na_X\psi+\frac{1}{n}\gamma(X)\D\psi.
\end{eqnarray}
\end{proposition}

The proof of Theorem~\ref{mainth} is based on the conformal covariance
of the fundamental Dirac operator $\D$ of the manifold $\M$; we now 
summarize some classical facts about Dirac operators in a conformal class of
the Riemannian metric $g$. So consider a nowhere vanishing function $h$ on
the manifold $\M$, and let $\bg=h^2 g$ be a conformal change of the metric. Then we have
an obvious identification between the two $\SO_n$-principal bundles of $g$
and $\bg$-orthonormal oriented frames denoted respectively by $\SO(\M)$ and
$\SO(\overline{\M})$. We can thus identify the corresponding $\Spin_n$-principal bundles
$\Spin(\M)$ and $\Spin(\overline{\M})$ and this leads to a bundle
isometry
\begin{equation}
\begin{array}{ccc}\label{cc}
\Sigma\M & \longrightarrow & \Sigma\overline{\M} \\
\varphi & \longmapsto & \overline{\varphi}.
\end{array}
\end{equation}

For more details, we refer to \cite{hitchin:74}, \cite{hijazi:86} or \cite{h1}. We can also relate the corresponding
Levi-Civita connections and Clifford multiplications. Indeed, denoting by $\overline{\na}$ and
$\overline{\gamma}$ the associated data which act on the bundle
$\Sigma\overline{\M}$. We can easily show that 
\begin{eqnarray}\label{cclcc}
\overline{\gamma}=h\gamma,\qquad\overline{\na}_X\psi-\na_X\psi=-\frac{1}{2h}\gamma(X)\gamma(\na
h)\psi-\frac{1}{2h}g(\na h,X)\psi.
\end{eqnarray}  

A result due to Hitchin \cite{hitchin:74} gives the conformal covariance of
the Dirac operator:

\begin{proposition}\label{prop5}
Let $\D$ (resp. $\overline{\D}$) be the fundamental
Dirac operator on the manifold $(\M^n,g)$ (resp. $(\M^n,\bg)$), then we
have the following identity:
\begin{eqnarray}
\overline{\D}(\overline{\psi})=h^{-\frac{n+1}{2}}\;\overline{\D(h^{\frac{n-1}{2}}\psi)},
\end{eqnarray}

for all $\psi\in\Ga(\Sigma\M)$.
\end{proposition} 

It is an obvious fact that this conformal change of the metric induces a
conformal change of the metric on the boundary. We can then identify the
connections and the Clifford multiplications of $\mathbf{S}(\pa\M)$ and
$\mathbf{S}(\overline{\pa\M})$. In the same way, the boundary Dirac operators
$\D^{\S}$ and $\DSC$, acting respectively on $\S(\pa\M)$ and $\S(\ov{\pa\M})$, satisfy
\begin{eqnarray}\label{p5}
\DSC\,(\overline{\psi})=h^{-\frac{n}{2}}\;\overline{\DS(h^{\frac{n-2}{2}}\psi)},
\end{eqnarray}

\noindent for all $\psi\in\Gamma\big(\S(\pa\M)\big)$. For more details on
these identifications, we refer to \cite{hijazi.montiel.zhang:2}.


\section{Local elliptic boundary conditions for the Dirac operator}\label{lebc}


In order to prove Theorem~\ref{mainth}, we have to use suitable boundary conditions for the fundamental Dirac operator
$\D$ on the manifold $\M$. In other words, we look for conditions
\begin{eqnarray*}
\mathbb{B}: \L^2(\mathbf{S}(\pa\M)) & \longrightarrow & \L^2(\V),
\end{eqnarray*}

where $\V$ is a Hermitian vector bundle over the boundary $\pa\M$, 
to impose on the restrictions of 
spinor fields on $\M$ to the boundary $\pa\M$ such that the Dirac operator
is a Fredholm operator, i.e. for given data $\Phi\in\Ga(\Sigma\M)$ and $\chi\in\Ga(\V)$ the following boundary value problem 
$$\left\lbrace
\begin{array}{ll}
\D\psi=\Phi & \quad\text{on}\;\M\\
\mathbb{B}\;(\psi_{|\pa\M})=\chi & \quad\text{along}\;\pa\M,
\end{array}
\right.\quad (BP)$$

has a unique solution up to a finite dimensional kernel. Moreover the
following eigenvalue problem 
$$\left\lbrace
\begin{array}{ll}
\D\varphi=\lambda\varphi & \quad\text{on}\;\M\\
\mathbb{B}\;(\varphi_{|\pa\M})=0 & \quad\text{along}\;\pa\M
\end{array}
\right.\quad(EBP)$$ 

should have a discrete spectrum with finite dimensional eigenspaces, unless
it is the whole complex plane. The preceding properties are satisfied if
the operator $\mathbb{B}$ satisfies some ellipticity conditions. We follow
\cite{hijazi.montiel.roldan:01} for the notion of ellipticity of a boundary
condition for the Dirac operator (for more general cases, we refer to
\cite{hormander} and \cite{lo}). In fact, the principal tool for finding well-posed
conditions was discovered by Cald\'eron and is called the Cald\'eron
projector of the Dirac operator $\D$, denoted by $\mathcal{P}_+(\D)$. This
projector is a pseudo-differential operator of order zero which has the
particularity that its principal symbol $\sigma\big(\mathcal{P}_+(\D)\big)$ detects ellipticity, i.e. such that problems $(BP)$ and
$(EBP)$ could be solved. Indeed, in \cite{bw:93}, the authors show that a
pseudo-differential operator 
\begin{eqnarray*}
\mathbb{B}:\L^2(\Sigma\M_{|\pa\M})\longrightarrow\L^2(\mathrm{V})
\end{eqnarray*}

defines an elliptic boundary condition for the operator $\D$ if it satisfies
the following property:
\begin{eqnarray*}
\sigma(\mathbb{B})(u)_{|\mathrm{Im}\;\sigma(\mathcal{P}_+(\D))(u)}:\mathrm{Im}\;\sigma(\mathcal{P}_+(\D))(u)\subset\Sigma_p\M\longrightarrow\mathrm{V}_p
\end{eqnarray*}

is an isomorphism on its image for all nontrivial $u\in\mathrm{T}_{p}(\pa\M)$
and all $p\in\pa\M$. Moreover if the rank of the vector bundle $\V$
and the dimension of $\mathrm{Im}\;\sigma(\mathcal{P}_+(\D))$ are the same 
then the boundary condition is said to be local. As the principal symbol of
the Cald\'eron projector of the Dirac operator is given by (see \cite{bw:93})
\begin{eqnarray*}
\sigma(\mathcal{P}_+(\D))(u)=\frac{1}{2|u|}(i\gamma^{\S}(u)+|u|\Id)
\end{eqnarray*}

for each nontrivial $u\in\mathrm{T}_{p}(\pa\M)$ and $p\in\pa\M$, we obtain
the following result which gives ellipticity of a boundary condition for
the Dirac operator $\D$ on a manifold with boundary (see \cite{hijazi.montiel.roldan:01}):

\begin{proposition}\label{p7}
Let $(\M^n,g)$ be an $n$-dimensional compact Riemannian spin manifold with
non-empty boundary $\pa\M$. Then, a pseudo-differential operator 
\begin{eqnarray*}
\mathbb{B}:\L^2(\mathbf{S}(\pa\M))\longrightarrow\L^2(\V),
\end{eqnarray*}

where $\V\rightarrow\pa\M$ is a Hermitian vector bundle, is an elliptic
boundary condition for the Dirac operator $\D$ of
$\M$ if and only if its principal symbol 
\begin{eqnarray*}
\sigma(\mathbb{B}):\TBM\longrightarrow\Hom_{\mathbb{C}}(\mathbf{S}(\pa\M),\mathrm{V})
\end{eqnarray*}

satisfies the following two conditions
\begin{enumerate}

\item $\ker\, \sigma(\mathbb{B})(u)\cap\{\varphi\in\Sigma_p\M\;/\;i\gamma(\nu)\gamma(u)\varphi=-|u|\varphi\}=\{0\},$

\item $\dim\,\im\, \sigma(\mathbb{B})(u)=\frac{1}{2}\dim\,\Sigma_p\M=2^{[\frac{n}{2}]-1}.$
\end{enumerate}

Moreover, if $\V$ is a bundle with rank
$\frac{1}{2}\dim\,\Sigma_p\M=2^{[\frac{n}{2}]-1}$, we have a local elliptic
boundary condition. When these ellipticity conditions are satisfied, the
following eigenvalue boundary problem 
$$\left\lbrace
\begin{array}{ll}
\D\psi=\lambda\psi & \quad\mathrm{on}\;\M\\
\mathbb{B}(\psi_{|\pa\M})=0 & \quad\mathrm{along}\;\pa\M,
\end{array}
\right.$$

has a discrete spectrum with
finite dimensional eigenspaces consisting of smooth spinor fields, unless
it is the whole complex plane.
\end{proposition}

We are now ready to study elliptic boundary conditions for the Dirac operator.


\section{The condition associated with a chirality operator}\label{chibound}


In this section, we consider an $n$-dimensional compact Riemannian spin
manifold $(\M^n,g)$ with non-empty boundary equipped with a {\it chirality
  operator}. First recall the definition of such an operator. A linear map
\begin{eqnarray*}
\G:\Ga(\Si\M)\longrightarrow\Ga(\Si\M)
\end{eqnarray*}

is a chirality operator if it satisfies the following properties:
\begin{equation}\label{poc}
\begin{array}{rc}
\G^2=\Id, & \<\G\varphi,\G\psi\>=\<\varphi,\psi\>\\
\na_X(\G\psi)=\G(\na_X\psi), & \gamma(X)\G\psi=-\G(\gamma(X)\psi)
\end{array}
\end{equation}

for all $\varphi$, $\psi\in\Ga(\Si\M)$ and $X\in\Ga(\TM)$. Note that such an
operator does not exist on all manifolds. However, we can note that if the
dimension $n=2m$ of the manifold $\M$ is even, then $\G=\gamma(\omega_{2m})$,
where $\omega_{2m}$ is the volume element of the spinor bundle, is a chirality operator. Now consider the fiber preserving
endomorphism
\begin{eqnarray*}
\gamma(\nu)\G:\Ga(\S(\pa\M))\longrightarrow\Ga(\S(\pa\M))
\end{eqnarray*}

acting on the restriction of the spinor bundle $\Si\M$ to the boundary. We
can easily check that this map is pointwise self-adjoint and is involutive.
So the bundle $\S(\pa\M)$ can be decomposed into two eigensubbundles $\V^{\pm}$
associated with the eigenvalues $\pm 1$. Now we can define
the two pointwise orthogonal projections:
$$\begin{array}{lccl}
\BCHI^{\pm}: & \L^2(\mathbf{S}(\pa\M)) & \longrightarrow & \L^2(\V^{\pm})\\
 & \varphi & \longmapsto & \frac{1}{2}(\Id\pm \gamma(\nu)\G)\varphi,
\end{array}$$ 

and it is easy to check that these operators satisfy the ellipticity
conditions given in Proposition~\ref{p7} (see \cite{hijazi.montiel.roldan:01}). 

\begin{remark}\label{rchi}
{\rm It is an important fact to note that a chirality operator $\G$ acting on
sections of $\Si\M$ allows to construct a chirality operator $\overline{\G}$
acting on $\Si\overline{\M}$. Indeed, the operator defined by
$$\begin{array}{clll}\label{cco}
\overline{\G}: & \Si\overline{\M} & \longrightarrow & \Si\overline{\M}\\
& \overline{\psi} & \longmapsto & \overline{\G}(\overline{\psi}):=\overline{\G\psi}
\end{array}$$

is a chirality operator acting on $\Si\overline{\M}$. }
\end{remark}

We can now prove the first part of the Theorem~\ref{mainth}:

\begin{theorem}\label{hijchi}
Let $(\M^n,g)$ be an $n$-dimensional compact Riemannian spin manifold with
non-empty boundary $\pa\M$ and with positive Yamabe invariant. Under the $\CHI$ boundary condition, the
spectrum of the Dirac operator $\D$ of $\M$ is a non-decreasing sequence of
real numbers $\{\lambda_k^{\CHI}\;/\;k\in\mathbb{Z}\}$ which satisfies
\begin{eqnarray*}
(\lambda^{\CHI}_k)^2\geq\frac{n}{4(n-1)}\,\mu_1(\L).
\end{eqnarray*}

\noindent Moreover, equality holds if and only if $\M$ is conformally
equivalent to the half-sphere $\hs(r)$, where $r$ depends of the first
eigenvalue of $\D$.
\end{theorem}

\pf
The spectrum is real because under this boundary condition, the Dirac
operator is symmetric. Indeed, if $\varphi$ and $\psi$ satisfy
$\BCHI^{\pm}(\varphi_{|\pa\M})=\BCHI^{\pm}(\psi_{|\pa\M})=0$, then we have:
\begin{eqnarray*}
\<\gamma(\nu)\psi,\varphi\>=\<\G(\gamma(\nu)\psi),\G\varphi\>=-\<\gamma(\nu)\psi,\varphi\>,
\end{eqnarray*}

hence by Formula~(\ref{ipp}), the symmetry property follows by integration. In fact, we can show that under this boundary condition, the Dirac operator extends to a self-adjoint linear operator on $\mathrm{L}^2$ (see \cite{wd} for this boundary condition or \cite{hormander} for a more general case). Furthermore, we have seen in
Section~\ref{lebc} that the eigenvalue boundary problem
\begin{equation}
\left\lbrace
\begin{array}{ll}\label{prochi}
\D\psi=\lambda^{\CHI}_k\psi & \text{on}\;\M\\
\BCHI^{\pm}(\psi_{|\pa\M})=0 & \text{along}\;\pa\M,
\end{array}
\right.
\end{equation}

admits a smooth solution $\psi\in\Ga(\Si\M)$ because of the ellipticity
of the $\CHI$ boundary condition. Let $\ov{g}=f^{\frac{4}{n-2}}g$ be a
conformal change of the metric and consider the spinor field
$\varphi=f^{-\frac{n-1}{n-2}}\psi\in\Ga(\Si\M)$.  Using the
conformal covariance of the Dirac operator given in
Proposition~\ref{prop5}, the spinor field $\ov{\varphi}\in\Ga(\Si\ov{\M})$ satisfies
\begin{eqnarray}\label{vpc}
\overline{\D}(\overline{\varphi})=\lambda^{\CHI}_k f^{-\frac{2}{n-2}}\,\overline{\varphi}
\end{eqnarray} 

Now putting this spinor field in the spinorial Reilly inequality (\ref{ir})
expressed in the metric $\ov{g}$ gives
\begin{align}\label{irchi}
\int_{\M}\Big(\frac{1}{4}\ov{\R}|\ov{\varphi}|^2-\frac{n-1}{n}|\lambda^{\CHI}_k|^2
f^{-\frac{4}{n-2}} & |\ov{\varphi}|^2\Big) dv(\bar{g}) \\
& \leq \int_{\pa\M}\Big(\<\ov{\D}^{\S}(\ov{\varphi}),\ov{\varphi}\>-\frac{n-1}{2}\ov{\H}|\ov{\varphi}|^2\Big)ds(\bg)\nonumber,
\end{align}

where $\overline{\R}$ (resp. $\overline{\H}$) is the scalar
curvature (resp. the mean curvature) of the manifold $(\M^n,\bg)$ (resp. of
the boundary $(\pa\M,\bg)$). However for $n\geq 3$, we can express the
corresponding curvatures using the conformal Laplacian
\begin{eqnarray*}
\L u=\frac{4(n-1)}{n-2}\Delta u +\R u,
\end{eqnarray*}

and the conformal mean curvature operator
\begin{eqnarray*}
\B u=\frac{2}{n-2}\frac{\partial u}{\partial\nu}+\H u.
\end{eqnarray*}

Indeed, we have  
\begin{eqnarray}\label{csc}
\overline{\R}=f^{-\frac{n+2}{n-2}}\;\L f,\qquad
\overline{\H}=f^{-\frac{n}{n-2}}\;\B f, 
\end{eqnarray}

where the function $f$ is the conformal factor of the metric $\bg$. 
Now recall that the eigenvalue boundary problem 
$$\left\lbrace
\begin{array}{ll}
\L u=\mu_1(\L)u &\quad\text{on } \M\\
\B u=0 &\quad\text{along }\;\pa\M,
\end{array}
\right.$$  

appearing in the statement of this theorem, was introduced by Escobar in
\cite{escobar:92} in the context of the Yamabe problem for manifolds with
boundary. He proved that the sign of the corresponding first eigenvalue
$\mu_1(\L)$, whose variational characterization is given by 
\begin{eqnarray*}
\mu_1(\L)= \underset{u\in C^{1}(\M),u\neq
  0}{\mathrm{inf}}\frac{\int_{\M}\big(\frac{2}{n-2}|\na
  u|^2+\frac{1}{2(n-1)}\R u^2\big)dv(g)+\int_{\pa\M}\H u^2
  ds(g)}{\int_{\M}u^2 ds(g)},
\end{eqnarray*}

is invariant under conformal change of the metric on $\M$ and that an
associated eigenfunction $f$ has to be positive. Moreover, Escobar showed
that $\mu_1(\L)$ has to be positive (resp. zero or negative) if and only if
there exists a conformally related metric on $\M$ with positive (resp. zero
or negative) scalar curvature and such that the boundary is minimal. Now
choosing the conformal factor of $\ov{g}$ to be a positive eigenfunction
$f_1$ associated with $\mu_1(\L)$ in Inequality (\ref{irchi}) and using the
relations (\ref{csc}) lead to
\begin{eqnarray}\label{chi25}
\int_{\M}\Big(\frac{1}{4}\mu_1(\L)-\frac{n-1}{n}|\lambda^{\CHI}|^2\Big)
|\overline{\varphi}|^2 f_1^{-\frac{4}{n-2}}dv(\bg)\leq\int_{\pa\M}\<\overline{\D}^{\S}(\overline{\varphi}),\overline{\varphi}\>ds(\bg).
\end{eqnarray}

We now prove that the boundary term vanishes under the $\CHI$ boundary
condition. Indeed, using the conformal covariance of the boundary Dirac operator
given in (\ref{p5}), we have
\begin{eqnarray*}
\overline{\D}^{\S}(\overline{\varphi})=f^{-\frac{n}{n-2}}_1\;\overline{\DS(f_1^{-\frac{1}{n-2}}\psi)}.
\end{eqnarray*}

Note that the volume forms of the boundary $\pa\M$ in the metric $g$ and
$\bg=f_1^{\frac{4}{n-2}}g$ are related by the formula
\begin{eqnarray*}
ds(\bg)=f_1^{\frac{2(n-1)}{n-2}}ds(g),
\end{eqnarray*} 

and then the boundary term is given by
\begin{eqnarray*}
\int_{\pa\M}\<\overline{\D}^{\S}(\overline{\varphi}),\overline{\varphi}\>ds(\bg)
& = &
\int_{\pa\M}\<\DS(f_1^{-\frac{1}{n-2}}\psi),f_1^{-\frac{1}{n-2}}\psi\>ds(g)\\
& = &
\int_{\pa\M}\<\gamma(d(f_1^{-\frac{1}{n-2}}))\psi,f_1^{-\frac{1}{n-2}}\psi\>ds(g)+\int_{\pa\M}
f_1^{-\frac{2}{n-2}}\<\DS\psi,\psi\>ds(g).
\end{eqnarray*}

We pointed out that the first term of the preceding identity is purely imaginary,
so inequality~(\ref{chi25}) gives
\begin{eqnarray}
\int_{\M}\Big(\frac{1}{4}\mu_1(\L)-\frac{n-1}{n}|\lambda^{\CHI}_k|^2\Big)
|\overline{\varphi}|^2 f_1^{-\frac{4}{n-2}}dv(\bg)\leq \int_{\pa\M}f_1^{-\frac{2}{n-2}}\<\DS\psi,\psi\>ds(g).
\end{eqnarray}

Using the fact that the spinor $\psi$ satisfies the eigenvalue boundary problem (\ref{prochi})
and that $\G\DS=\DS\G$ (using property (\ref{poc})), we obtain:
\begin{eqnarray*}
\<\DS\psi,\psi\>=\<\gamma(\nu)\G(\DS\psi),\gamma(\nu)\G\psi\>=-\<\DS\psi,\psi\>,
\end{eqnarray*}

hence, the desired inequality follows. Assume now that equality is
achieved. So equality occurs in (\ref{irchi}) and then by
Proposition~\ref{p4}, 
the spinor field $\ov{\varphi}$ satisfies the following equation
\begin{eqnarray*}
\ov{\Pe}_X(\ov{\varphi})=\ov{\na}_X\ov{\varphi}+\frac{1}{n}\ov{\gamma}(X)\ov{\D}(\ov{\varphi})=0,\qquad\forall
X\in\Gamma(\TM).
\end{eqnarray*}

Moreover, using equation~(\ref{vpc}), we conclude that this spinor field satisfies the generalized Killing equation
\begin{eqnarray*}
\ov{\na}_X\ov{\varphi}+\frac{\lambda}{n}f_1^{-\frac{2}{n-2}}\;\ov{\gamma}(X)\ov{\varphi}=0,
\end{eqnarray*}

for all $X\in\Gamma(\TM)$ and where $\frac{\lambda}{n}f_1^{-\frac{2}{n-2}}$
is a real-valued function. However, it is a well-known fact (see
\cite{hijazi:86}) that, because $f_1$ is a real-valued function, $f_1$ is constant. Then the
spinor field $\varphi$ (and so the spinor field $\psi$) is a real Killing
spinor, i.e. it satisfies the Killing equation
\begin{eqnarray*}
\na_X\psi+\frac{c}{2}\gamma(X)\psi=0,\qquad\forall
X\in\Gamma(\TM),
\end{eqnarray*}

where $c$ is a positive real number given by
$c=\frac{2\lambda}{n}$. Moreover this implies that the length
$|\psi|^2$ is a non-zero constant and that $(\M,g)$ is an Einstein
manifold with Ricci curvature $\mathrm{Ric}=(n-1)c^2 g$ whose boundary is
minimal. Now consider the function given by
$\F=\<\G(\psi),\psi\>$ which is real-valued because the chirality
operator $\G$ is pointwise self-adjoint. We then check that $\F$ is non
identically zero on $\M$ since using Formula~(\ref{ipp}), we have
\begin{eqnarray*}
nc\int_{\M}\F dv(\ov{g})=\int_{\pa\M}|\psi|^2ds(\ov{g}),
\end{eqnarray*}

and then, as $\psi$ is a non-zero constant lenght spinor field, we
obtain $c>0$ and $\F\not\equiv 0$. We now prove that the function $\F$ satisfies the boundary problem
$$\left\lbrace
\begin{array}{ll}
\Delta\F=nc^2\F & \quad\text{on}\;\M\\
\F_{|\pa\M}=0 & \quad\text{along}\;\pa\M.
\end{array}
\right.$$

An easy calculation using the Killing equation gives $\Delta\F=nc^2\F$ on
$\M$. The spinor field $\psi$ satisfies the eigenvalue boundary
problem~(\ref{prochi}), in particular we have
\begin{eqnarray*}
\ga(\nu)\G\psi_{|\pa\M}=\mp\psi_{|\pa\M}
\end{eqnarray*}
 
and then
\begin{eqnarray*}
\F_{|\pa\M}=\<\G\psi_{|\pa\M},\psi_{|\pa\M}\>=\pm\<\ga(\nu)\psi_{|\pa\M},\psi_{|\pa\M}\>,
\end{eqnarray*}

hence $\F_{|\pa\M}=0$ since the right hand side of the last equation is
purely imaginary. Applying the Obata Theorem for manifold with boundary
proved by Reilly in \cite{reilly} allows to conclude that $\M$ is conformally
equivalent to $\hs(r)$ with $r=\frac{1}{c}$. 
\qed\\

The proof of Theorem~\ref{hijchi} is based on a well-chosen conformal
metric $\ov{g}=f^{\frac{4}{n-2}}g$ which has no sense if $n=2$. However, we
can apply with slight modifications, the argument used in \cite{baer:92}
for the case of compact surfaces without boundary.

\begin{theorem}\label{hijchi2}
Let $(\M^2,g)$ be a compact Riemannian spin surface of genus $0$. Suppose
$\pa\M$ is compact connected. Under the $\CHI$ boundary condition, the
spectrum of the Dirac operator $\D$ of $\M$ is a non-decreasing sequence of
real numbers $\{\lambda_k^{\CHI}\;/\;k\in\mathbb{Z}\}$ which satisfies
\begin{eqnarray}\label{surchi}
\Big(\lambda^{\CHI}_k\Big)^2\geq\frac{2\pi}{\aire}. 
\end{eqnarray}

\noindent Equality holds for the smallest eigenvalue if and only if $\M^2$ is
isometric to a standard hemisphere $\mathbb{S}^2_+(r)$ with
$r=\frac{1}{\lambda^{\CHI}_{\pm 1}}$.
\end{theorem}

\pf
We first show that under this boundary condition, any eigenvalue
$\lambda^{\CHI}$ of the Dirac operator satisfies 
\begin{eqnarray}\label{sur}
\big(\lambda^{\CHI}\big)^2\geq\frac{1}{2}\sup_{u}\inf_{\M}\left(\ov{\R}\,e^{2u}\right),
\end{eqnarray}

\noindent where the supremum is taken over all the functions $u$ satisfying
$\frac{\pa u}{\pa\nu}+\H=0$, where $\R$ (resp. $\ov{R}$) is the scalar
curvature of $g$ (resp. $\ov{g}=e^{2u}g$) and $\H$ is the
geodesic curvature of $\pa\M$. So for a metric $\ov{g}=e^{2u}g$ with $\frac{\pa u}{\pa\nu}+\H=0$, in the
conformal class of $g$, the associated Dirac operators are related by the
following identity
\begin{eqnarray*}
\ov{\D}(e^{-\frac{1}{2}u}\ov{\psi})=e^{-\frac{3}{2}u}\ov{\D\psi}.
\end{eqnarray*}  
Using the argument given in the proof of Theorem~\ref{hijchi} and the
fact that $\ov{H}=0$ for all $u$ such that $\frac{\pa u}{\pa\nu}+\H=0$ lead to
the inequality
\begin{eqnarray*}
\int_{\M}\left(\frac{\ov{R}}{4}e^{2u}-\frac{|\lambda^{\MIT}_k|^2}{2}\right)e^{-2u}|\ov{\varphi}|^2dv(\ov{g})\leq\int_{\pa\M}\<\ov{\DS}(\ov\varphi),\ov{\varphi}\>ds(\ov{g}),
\end{eqnarray*}

where $\varphi=e^{-\frac{1}{2}u}\psi$ and $\psi\in\Gamma(\Sigma\M)$
satisfies the boundary problem (\ref{prochi}). The conformal covariance of
the $\CHI$ boundary condition gives (\ref{sur}). We now give an explicit
calculation of the right-hand side of Inequality (\ref{sur}). First, remark that for
$n=2$, the transformations of the scalar curvature and the geodesic
curvature give
$$\left\lbrace
\begin{array}{lll}
\ov{\R}e^{2u} & = & \R+2\Delta u\\
\ov{H}e^u & = & \frac{\pa u}{\pa\nu}+\H,
\end{array}
\right.$$

and then, by assumption on the function $u$, $\ov{\H}=0$. Now, note that for
such a function $u$
\begin{eqnarray*}
\inf_{\M}(\ov{\R}e^{2u}) & \leq &
\frac{1}{\aire}\Big(\int_{\M}\ov{\R}e^{2u}dv(g)+2\int_{\pa\M}\ov{\H}ds(g)\Big)\\
& \leq & \frac{1}{\aire}\Big(\int_{\M}\R dv(g)+2\int_{\pa\M}\H ds(g)\Big)=\frac{4\pi\chi(\M^2)}{\aire}
\end{eqnarray*}

by Stokes and Gauss-Bonnet Formula for surfaces with boundary (see
\cite{dc}). Let $u_1$ a solution of the boundary problem (see \cite{t} for example)
\begin{equation}
\left\lbrace
\begin{array}{llll}
2\Delta u_1 = \frac{1}{\aire}\Big(\int_{\M}\R dv(g)+2\int_{\pa\M}\H
ds(g)\Big)-\R & \quad\text{on}\;\M\nonumber\\
\frac{\partial u_1}{\partial\nu}+\H = 0 & \quad\text{along}\;\pa\M,
\end{array}
\right.
\end{equation}

then an easy calculation gives
$$\inf_{\M}(\ov{\R}e^{2u_1})=\frac{4\pi\chi(\M^2)}{\aire}.$$

From Inequality~(\ref{sur}), we obtain:
\begin{eqnarray*}
\Big(\lambda^{\CHI}_k\Big)^2\geq\frac{2\pi\chi(\M^2)}{\aire}.
\end{eqnarray*}

However, the surface $\M$ is of genus $0$ and its boundary has one
connected component, then $\chi(\M)=1$ and so Inequality (\ref{surchi})
follows immediately. The equality case is treated as in the proof of Theorem~\ref{hijchi} .
\qed\\

\begin{remark}
\rm{The Euler-Poincar\'e characteristic of a surface of genus $g\geq0$ and with
$m\geq 1$ boundary components is given by (see~\cite{hirsch})
\begin{eqnarray*}
\chi(\M)=2-2g-m.
\end{eqnarray*}

It is a simple fact that
\begin{eqnarray*}
\chi(\M)>0\Leftrightarrow\chi(\M)=1\Leftrightarrow\left(g=0\quad\text{et}\quad
m=1\right).
\end{eqnarray*}}
\end{remark}

We can now relate the Yamabe invariant of the manifold $\M$ with the
eigenvalues of the Dirac operator under the boundary condition associated
with a chirality operator. This conformal invariant, denoted by $\mu(\M)$, has been introduced by
Escobar in \cite{escobar:92} in order to solve the Yamabe problem for
manifolds with boundary and has the following variational characterization 

\begin{eqnarray*}
\mu(\M)=\underset{u\in C^{1}(\M),u\neq
  0}{\mathrm{inf}}\frac{\int_{\M}\big(\frac{2}{n-2}|\na
  u|^2+\frac{1}{2(n-1)}\R u^2\big)dv(g)+\int_{\pa\M}\H u^2
  ds(g)}{\Big(\int_{\M}u^{\frac{2n}{n-2}} ds(g)\Big)^{\frac{n-2}{n}}}.
\end{eqnarray*}

He proved that $\mu(\M)$ has the same sign as $\mu_1(\L)$ and it is
invariant with respect to conformal changes of the metric on $\M$. The
H\"older inequality applied to an eigenfunction $f_1$ associated with
$\mu_1(\M)$ gives
\begin{eqnarray*}
\mu_1(\L)\geq\frac{\mu(\M)}{\mathrm{vol}(\M,g)^{\frac{2}{n}}}
\end{eqnarray*}

and equality implies that $f_1$ is constant. Thus, from
Theorem~\ref{hijchi} and Theorem~\ref{hijchi2}, we have:
\begin{corollary}
Let $(\M^n,g)$ be an $n$-dimensional compact Riemannian spin manifold with
non-empty compact connected boundary $\pa\M$ and $n\geq 2$. Under the $\CHI$ boundary
condition, any eigenvalue $\lambda^{\CHI}$ of the Dirac operator
$\D$ satisfies 
\begin{eqnarray}
|\lambda^{\CHI}|^2\,\mathrm{vol}(\M,g)^{\frac{2}{n}}\geq\frac{n}{4(n-1)}\mu(\M). 
\end{eqnarray}
\end{corollary}

\begin{remark}\label{rem2}
\rm{Theorem~\ref{hijchi} and Theorem~\ref{hijchi2} can also be proved for the
$\MIT$ bag boundary condition used in \cite{hijazi.montiel.roldan:01} and  \cite{hijazi.montiel.zhang:2}
for example. This condition is defined as being the orthogonal projection
on the eigensubbundles of the pointwise endomorphism $i\ga(\nu)$ acting on
$\mathbf{S}(\pa\M)$ associated with the eigenvalues $\pm 1$, i.e
$$\begin{array}{lccl}
\BMIT^{\pm}: & \L^2(\mathbf{S}(\pa\M)) & \longrightarrow & \L^2(\V^{\pm})\\
 & \varphi & \longmapsto & \frac{1}{2}(\Id\pm i\gamma(\nu))\varphi.
\end{array}$$ 

This differential operator satisfies the Lopatinsky-Shapiro ellipticity
conditions, and so it defines a local elliptic boundary condition for the 
Dirac operator of the manifold $\M$. We can then show (see
\cite{hijazi.montiel.roldan:01}) that under the $\BMIT^-$ (resp. $\BMIT^+$)
condition, the spectrum of the Dirac operator is a discrete set of complex
numbers with positive (resp. negative) imaginary part. Furthemore, we can
show that any eigenvalue $\lambda^{\MIT}$ under this boundary condition
satifies for $n\geq 3$
\begin{eqnarray}\label{hijmit}
|\lambda^{\MIT}|^2>\frac{n}{4(n-1)}\,\mu_1(\L),
\end{eqnarray}

and for $n=2$
\begin{eqnarray}\label{hijmit1}
|\lambda^{\MIT}|^2>\frac{2\pi}{\aire}. 
\end{eqnarray}

Equality cannot be achieved in (\ref{hijmit}) and (\ref{hijmit1})
 otherwise there should exist on the manifold
$\M$ a generalized Killing spinor with imaginary Killing function and then
the scalar curvature should be non-positive. However, if equality holds,
the scalar curvature of $\M$ is positive and so there is a contradiction.}
\end{remark}


\bibliographystyle{amsalpha}     
\bibliography{biblioart1}

\providecommand{\bysame}{\leavevmode\hbox to3em{\hrulefill}\thinspace}
\providecommand{\MR}{\relax\ifhmode\unskip\space\fi MR }
\providecommand{\MRhref}[2]{%
  \href{http://www.ams.org/mathscinet-getitem?mr=#1}{#2}
}
\providecommand{\href}[2]{#2}
\begin{thebibliography}{BBW93}

\bibitem[B{\"a}r92]{baer:92}
C.~B{\"a}r, \emph{Lower eigenvalue estimate for {D}irac operator}, Math. Ann.
  \textbf{293} (1992), 39--46.

\bibitem[B{\"a}r98]{baer:98}
\bysame, \emph{Extrinsic bounds of the {D}irac operator}, Ann. Glob. Anal.
  Geom. \textbf{16} (1998), 573--596.

\bibitem[Bau89]{h1}
H.~Baum, \emph{Complete {R}iemannian manifolds with imaginary {K}illing
  spinors}, Ann. Glob. Anal. Geom. \textbf{7} (1989), 205--226.

\bibitem[BBW93]{bw:93}
B.~Boo\ss-Bavnbek and K.~P. Wojciechowski, \emph{Elliptic boundary problems for
  the {D}irac operator}, Birkh\"auser, Basel, 1993.

\bibitem[DC91]{dc}
M.~{D}o {C}armo, \emph{Differential forms and applications}, Springer-Verlag,
  1991.

\bibitem[DW95]{wd}
S.~D\"urr and A.~Wipf, \emph{Gauge theories in a bag}, Nucl. Phys. B
  \textbf{443} (1995), 201--232.

\bibitem[Esc92]{escobar:92}
J.~F. Escobar, \emph{The {Y}amabe problem on manifolds with boundary}, J. Diff.
  Geom. \textbf{35} (1992), 21--84.

\bibitem[Fri80]{friedrich}
T.~Friedrich, \emph{Der erste {E}igenwert des {D}irac-{O}perators einer
  kompakten {R}iemannschen {M}annigfaltigkeit nicht negativer
  {S}kalarkr\"ummung}, Math. Nach. \textbf{97} (1980), 117--146.

\bibitem[Fri00]{fried}
\bysame, \emph{Dirac operators in {R}iemannian geometry}, vol.~25, Amer. Math.
  Soc. Graduate Studies in Math., 2000.

\bibitem[Hij86]{hijazi:86}
O.~Hijazi, \emph{A conformal lower bound for the smallest eigenvalue of the
  {D}irac operator and {K}illing spinors}, Commun. Math. Phys. \textbf{25}
  (1986), 151--162.

\bibitem[Hij91]{hijazi:91}
\bysame, \emph{Premi\`ere valeur propre de l'op\'erateur de {D}irac et nombre
  de {Y}amabe}, C. R. Acad. Sci. Paris \textbf{313} (1991), 865--868.

\bibitem[Hir76]{hirsch}
M.~W. Hirsch, \emph{Differential {T}opology}, New {Y}ork: {S}pringer-{V}erlag,
  1976.

\bibitem[Hit74]{hitchin:74}
N.~Hitchin, \emph{Harmonic spinors}, Adv. Math. \textbf{14} (1974), 1--55.

\bibitem[HMR02]{hijazi.montiel.roldan:01}
O.~Hijazi, S.~Montiel, and S.~Rold\'an, \emph{Eigenvalue boundary problems for
  the {D}irac operator}, Commun. Math. Phys. \textbf{231} (2002), 375--390.

\bibitem[HMZ01]{hijazi.montiel.zhang:1}
O.~Hijazi, S.~Montiel, and X.~Zhang, \emph{Dirac operator on embedded
  hypersurfaces}, Math. Res. Lett. \textbf{8} (2001), 20--36.

\bibitem[HMZ02]{hijazi.montiel.zhang:2}
\bysame, \emph{Conformal lower bounds for the {D}irac operator on embedded
  hypersurfaces}, Asian J. Math. \textbf{6} (2002), 23--36.

\bibitem[H{\"o}r85]{hormander}
L.~H{\"o}rmander, \emph{The analysis of linear partial differential operators
  {III}}, Springer, Berlin, 1985.

\bibitem[LM89]{lawson.michelson:89}
H.~B. Lawson and M.~L. Michelsohn, \emph{Spin {G}eometry}, {P}rinceton
  {U}niversity {P}ress ed., vol.~38, Princeton Math. Series, 1989.

\bibitem[Lop53]{lo}
Ya.~B. Lopatinski$\breve{\i}$, \emph{On a method for reducing boundary problems
  for a system of differential equations of elliptic type to regular integral
  equations}, Ukrain. Math. $\check{\mathrm{Z}}$. \textbf{5} (1953), 123--151,
  (Russian).

\bibitem[Mor01]{morel}
B.~Morel, \emph{Eigenvalue estimates for the {D}irac-{S}chr{\"o}dinger
  operators}, Journal of Geometry and Physics \textbf{38} (2001), 1--18.

\bibitem[Rei77]{reilly}
R.~C. Reilly, \emph{Application of the {H}essian operator in a {R}iemannian
  manifold}, Indiana Univ. Math. J. \textbf{26} (1977), 459--472.

\bibitem[Tay96]{t}
M.~Taylor, \emph{Partial differentiel equations, {V}ol. 1: {B}asic theory},
  Springer, 1996.

\bibitem[Tra95]{trautman}
A.~Trautman, \emph{The {D}irac operator on hypersurfaces}, Acta Phys. Polon.
  \textbf{6} (1995), 1283--1310.

\bibitem[Yam60]{yamabe:60}
H.~Yamabe, \emph{On a deformation of {R}iemannian structures on compact
  manifolds}, Osaka Math. J. (1960), 21--37.

\end{thebibliography}


\vspace{0.8cm}     
Author address:     
\nopagebreak     
\vspace{5mm}\\     
\parskip0ex     
\vtop{     
\hsize=6cm\noindent     
\obeylines     

}     
\vtop{     
\hsize=8cm\noindent     
\obeylines     
Simon Raulot,
Institut \'Elie Cartan BP 239     
Universit\'e de Nancy 1     
54506 Vand\oe uvre-l\`es -Nancy Cedex     
France     
}     
     
\vspace{0.5cm}     
     
E-Mail:     
{\tt raulot@iecn.u-nancy.fr  }

\end{document}